\documentclass[english,a4paper,10pt]{article}
\usepackage[utf8]{inputenc}
\usepackage[T1]{fontenc}
\usepackage{amsmath}
\usepackage{amsthm}
\usepackage{graphicx}
\usepackage{amssymb}
\usepackage{esint}
\usepackage{upref} \usepackage{hyperref}

\newcommand \Pl {\mathbb{P}l}

\newcommand \jij {J(i,j,\theta^2)}
\newcommand \poissdev {\gamma|\theta-\sqrt{n}|}
\newtheorem*{theoremnn}{Theorem}
\newtheorem{theorem}{Theorem}
\newtheorem{deplemma}{Depoissonization Lemma}
\newtheorem{proposition}{Proposition}[section]

\title{Variance of Linear Statistic for Plancherel Young Diagrams}
\date{}
\author{Konstantin Tolmachov}

\begin{document}
\maketitle
\begin{abstract}
In this paper we compute the precise asymptotics of the variance of linear
statistic of descents on a growing interval for Plancherel Young diagrams
(following Vershik and Kerov, diagrams are considered rotated by $\pi/4$). We
also give an example of a local configuration with linearly growing variance in
a fixed regime and prove the central limit theorem for this configuration in the
given regime.\end{abstract}
\section{Introduction}
\emph{Young diagram} $\lambda$ with $n$ cells is a table, rows of which
represent a partition of $n$ to the sum of nonincreasing summands
$\lambda_{1}\geq\lambda_{2}\geq...\geq\lambda_{m}>0,\,\lambda_{1}+...+\lambda_{m
}=n$.
Set $n=|\lambda|$. We denote $\mathbb{Y}_{n}$ the set of Young diagrams with $n$
cells. Following Vershik and Kerov we will draw diagrams rotated by
$\dfrac{\pi}{4}$:

\includegraphics{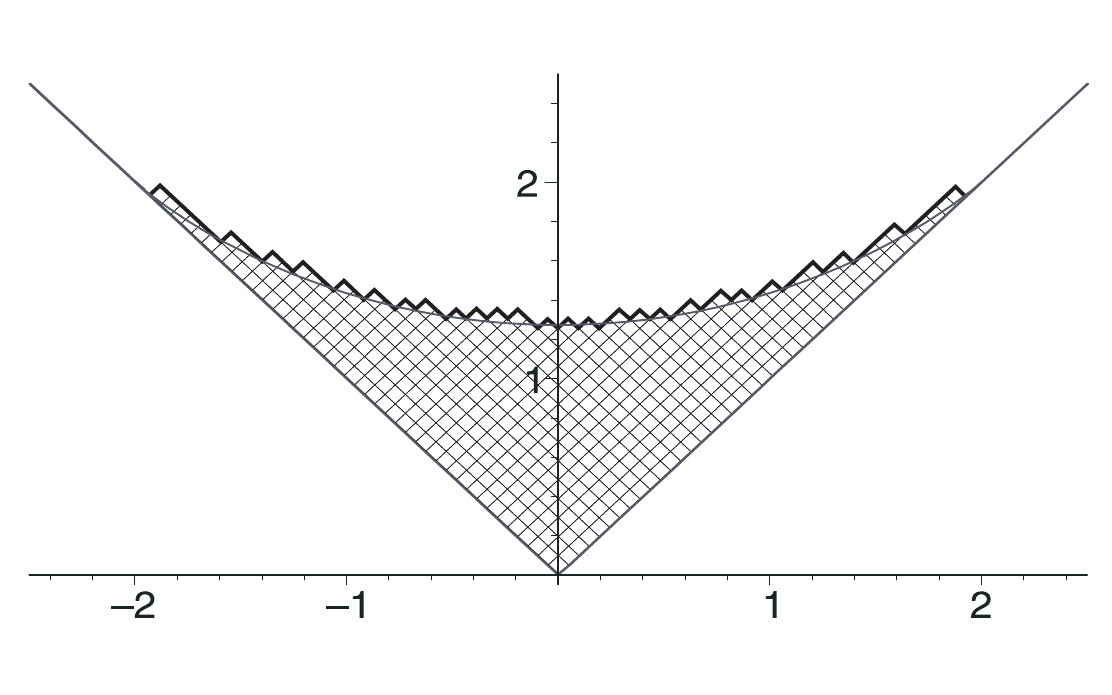}

For a diagram $\lambda$ define a sequence
$(c_{i}(\lambda))\in\{0,1\}^{\mathbb{Z}}$
as follows:
$$
\begin{cases}
c_{k}(\lambda)=1, & k=\lambda_i-i, \text{for some } i;\\
c_{k}(\lambda)=0, & \text{otherwise}.\end{cases}
$$

This sequence has the following geometrical interpretation. Consider a function
$\Phi_{\lambda}$ whose graph represents the upper edge of a rotated diagram. The
sequence $c_{i}(\lambda)$ is the sequence of descents of a rotated diagram, that
is, $\Phi_{\lambda}$ is a continuous piecewise-linear function with the following
property:
$$
\begin{cases}
 \Phi_{\lambda}'(x)=1, & c_{[x]}(\lambda)=0;\\
 \Phi_{\lambda}'(x)=-1, & \text{otherwise}.\end{cases}
$$

The probability measure $\Pl^{(n)}$ on $\mathbb{Y}_{n}$ is given by the formula
$$
\Pl^{(n)}(\lambda)=\frac{(\dim\lambda)^{2}}{n!},
$$
where $\dim\lambda$ is the dimension of the irreducible representation of the
symmetric group $\mathbb{S}_{n}$ corresponding to $\lambda$. This measure is
called the \emph{Plancherel measure}.

It is shown in \cite{key-2} and \cite{key-8} that the local distribution of
$c_{i}(\lambda)$ as ${|\lambda|\to\infty}$ is governed by the discrete
determinantal sine-process. For the reader's convenience we recall the statement
of Theorem 2 from \cite{key-2}:
for an arbitrary subset $\{d_{1},...,d_{s}\}\subset\mathbb{Z}$ and a sequence
$(x_n),\lim\limits_{n \to \infty}x_n/\sqrt{n} = u \in (-2,2)$ the following
holds:
$$
\lim\limits_{n \to
\infty}\mathbb{E}_{\mathbb{P}l^{(n)}}(c_{x_n+d_1}(\lambda)c_{x_n+d_2}
(\lambda)...c_{x_n+d_s}(\lambda))=\det\left[K_{\sin}(d_i-d_j,\arccos(u/2))\right
]_{1\leq i,j\leq s},
$$
where $K_{\sin}(x,\phi)$ is the \emph{sine kernel},
$$
K_{\sin}(x,\phi)=\frac{\sin(\phi(x-y))}{\pi(x-y)}.
$$
(for the more general statement see \cite{key-2}).

In this paper we are concerned with the \emph{linear statistic}
$\displaystyle{\sum_{i \in I}c_i}$ of descents for Plancherel Young diagrams. The
main result of this paper is the precise asymptotics for the variance of the
linear statistic of $c_i$: 
\begin{theorem}
\label{main}
 For all $a,b\in(-2,2)$ and all sequences ${x_n},{y_n}$ such that
$\lim\limits_{n\to\infty} {\dfrac{x_n}{\sqrt{n}}=a}$, $\lim\limits_{n\to\infty}
\dfrac{y_n}{\sqrt{n}}=b$ and $\lim\limits_{n\to\infty} (y_n-x_n)=+\infty$, the
following holds:
\begin{equation}
\label{Var}
 \lim\limits_{n\to\infty}\frac{Var_{\Pl^{(n)}}\left(\displaystyle{\sum\limits_{
i=x_n}^{y_n}} c_{i}(\lambda)\right)}{\log(y_n-x_n)}=\frac{1}{\pi^2}.
\end{equation}
\end{theorem}

\paragraph{Plan of the proof and organization of the paper}
To prove Theorem \ref{main} we use the theorem established, independently and
simultaneously, by Borodin, Okounkov and Olshanski \cite{key-2} and Johansson
\cite{key-8}, which
claims that the poissonization of the Plancherel measure is the discrete Bessel
determinantal point process. In Section \ref{sec_poiss} we analyze the
asymptotics of the poissonized variance in the left part of (\ref{Var}). In
Section \ref{sec_depoiss} switch to the asymptotics with respect to the
Plancherel measure (depoissonization) is described. Estimates of the Bessel
kernel that we need in Section \ref{sec_poiss} are proved in Section
\ref{sec_est}.

Using general properties of determinantal point processes it is easy to obtain
an upper bound for the variance of frequency of general local patterns on Young
diagram. In Section \ref{corners} we give an example of local patterns with
linearly growing variance and prove the central limit theorem for them.

\paragraph{Acknowledgements} I am deeply grateful to Alexander Bufetov for the
statement of the problem and for many helpful discussions. I am deeply grateful
to Vadim Gorin, Sevak Mkrtchyan and Grigori Olshanski for helpful discussions.

\section{Poissonization}
To prove Theorem {\ref{main}} we use the \emph{poissonization}(\cite{key-2}, \cite{key-8}).
Set $\mathbb{Y}=\bigcup\limits_{n=1}^{\infty}\mathbb{Y}_n$,
and for $\eta \in \mathbb{C}$ let
$$
\Pl_{\eta}(\lambda)=e^{-\eta}\sum_{k=0}^{\infty}\frac{\eta^{k}}{k!}\mathbb{P}l^{
(k)}(\lambda),\,\eta\in\mathbb{C},
$$
be \emph{the poissonization of Plancherel measure}
(set $\Pl^{(k)}(\lambda)=0$ when $|\lambda|\neq k$).
Borodin, Okounkov, Olshanski in \cite{key-2} and Johansson in \cite{key-8}
showed that $\Pl_\eta$ induces the \emph{determinantal point process} on the set of
Young diagrams:
denote $\eta=\theta^2$; for an arbitrary subset
$\{x_{1},...,x_{s}\}\subset\mathbb{Z}$ the following holds:
$$
\mathbb{E}_{\mathbb{P}l_{\theta^2}}(c_{x_1}(\lambda)c_{x_2}(\lambda)...c_{x_s}
(\lambda))=\det\left[J(x_{i},x_{j},\theta^{2})\right]_{1\leq i,j\leq s},
$$
where $J(x,y,\theta^{2})$ is the \emph{Bessel kernel} defined by
$$
J(x,y,\theta^{2})=\theta\frac{J_{x}(2\theta)J_{y+1}(2\theta)-J_{x+1}(2\theta)J_{
y}(2\theta)}{x-y}.
$$
(from now on we understand variance and expectation with respect to
$\Pl_{\theta^2}$ for $\theta\in\mathbb{C}$ formally).

It was shown in \cite{key-2} and \cite{key-8} that local patterns of a Young
diagram are governed by the sine-process (i.e. the Bessel kernel degenerates to
the sine kernel as $n \to \infty$).

There are continuous versions of sine and Bessel processes, see \cite{key-13},
\cite{key-11}, \cite{key-12}. The central limit theorem for the continuous
sine-process was proved in \cite{key-11}. The central limit theorem for the
continuous Bessel process was proved in \cite{key-12}. Note that in both cases
one has a logarithmic growth of the variance of the linear statistics. See
\cite{key-13} for the central limit theorem and other general results for
determinantal point processes. 

Following \cite{key-1} and using the well-known identity 
$$
J(x,x,\theta^2)=\sum\limits_{y\in\mathbb{Z}}\left(J(x,y,\theta^2)\right)^2, 
$$
which represents the fact that the linear operator defined by the kernel $J$ is a projection (see \cite{key-2}), we
rewrite the poissonized variance from \ref{main} in the following way:
$$
Var_{\theta^2}\left(\sum\limits_{x}^{y}c_{i}(\lambda)\right) =
\sum\nolimits_{i\in[x_n,y_n]}\sum\nolimits_{j\notin[x_n,y_n]}
\left(\jij\right)^2,
$$
where $Var_{\theta^2}$ is the variance with respect to the poissonization of the
Plancherel measure $\Pl_{\theta^2}$.
Now we formulate the poissonized version of Theorem \ref{main}.
\begin{proposition}
 \label{poissmain}
There is a constant $\gamma>0$ such that for all $a,b\in(-2,2)$ and all
sequences ${x_n},{y_n}$ such that $\lim\limits_{n\to\infty}
\dfrac{x_n}{\sqrt{n}}=a$, $\lim\limits_{n\to\infty} \dfrac{y_n}{\sqrt{n}}=b$ and
${\lim\limits_{n\to\infty} (y_n-x_n)=\infty}$ the following holds:
\begin{equation}
\label{poissVar}
 \left|\frac{\displaystyle{\sum\nolimits_{i\in[x_n,y_n]}\sum\nolimits_{j\notin[
x_n,y_n]}}
\left(\jij\right)^2}{\log(y_n-x_n)}-\frac{1}{\pi^2}\right|\exp(-\gamma|\sqrt{n}
-\theta|)=o(1).
\end{equation}
\end{proposition}
Similar proposition is stated in \cite{key-5}.
Note that Proposition \ref{poissVar} immediately yields an asymptotical formula
for the variance of poissonized measure $Pl_n, {n > 0}$, and from the Soshnikov's
central limit theorem from  \cite{key-13} we get
\begin{proposition}
 For all $a,b\in(-2,2)$ and all sequences ${x_n},{y_n}$ such that
${\lim\limits_{n\to\infty} \dfrac{x_n}{\sqrt{n}}=a}$, ${\lim\limits_{n\to\infty}
\dfrac{y_n}{\sqrt{n}}=b}$ and $\lim\limits_{n\to\infty} (y_n-x_n)=+\infty$ the
following holds:
 $$
  \frac{\sum_{i \in [x_n,y_n]} c_i(\lambda) - \mathbb{E}_n (\sum_{i \in
[x_n,y_n]} c_i(\lambda))}{\frac{1}{\pi}\sqrt{\log (y_n-x_n)}} \xrightarrow[d]{} \mathcal{N}(0,1)
 $$
 as $n\to\infty$.
\end{proposition}
Here $\xrightarrow[d]{}$ means convergence in distribution.

The central limit theorem for the integral of the deviation of a Young diagram
from its limit shape is proved in \cite{key-9}, \cite{key-10}. The central limit
theorem for the point-wise deviation is stated in \cite{key-5}.

It is clear from the sine kernel asymptotics that expectation of linear
statistic $\mathbb{E}\left(\sum c_i\right)$ grows linearly in ${(y_n-x_n)}$.   

One may come back to the Plancherel measure from its poissonization using the
depoissonization lemma (Lemma 3.1 in \cite{key-2}) and some of its modifications
(Lemmas 3.1
- 3.2 in \cite{key-1}). 

Fix $0<\alpha<1/4$.

\begin{deplemma}
Let $\{f_{n}\}$ be a sequence of entire functions,
$$
f_{n}=e^{-z}\sum_{k=0}^{\infty}\frac{f_{nk}}{k!}z^{k},
$$
and assume that there exist constants $f_{\infty}$ and $\gamma$ such that 
$$
\max_{|z|=n}f_{n}(z)=O(e^{\gamma\sqrt{n}})
$$
and
$$
\max_{|z/n-1|\leq
n^{-\alpha}}|f_{n}(z)-f_{\infty}|e^{-\gamma|z-n|/\sqrt{n}}=o(1)
$$
 as $n\rightarrow\infty$. Then
$$
\lim_{n\rightarrow\infty}f_{nn}=f_{\infty}.
$$
\end{deplemma}

\begin{deplemma}
Let $\{f_{n}\}$ be a sequence of entire functions,
$$
f_{n}=e^{-z}\sum_{k=0}^{\infty}\frac{f_{nk}}{k!}z^{k},
$$
and assume that there exist constants $f_{\infty}$, $\gamma$, $C_{1}$,
$C_{2}$ such that
$$
\max_{|z|=n}|f_{n}(z)|\leq C_{1}e^{\gamma\sqrt{n}},
$$
$$
\max_{|z/n-1|\leq n^{-\alpha}}|f_{n}(z)-f_{\infty}|e^{-\gamma|z-n|/\sqrt{n}}\leq
C_{2}.
$$
Then there exists a constant $C=C(\gamma,C_{1},C_{2})$ such that
for all ${n>0}$ we have
$$
|f_{nn}-f_{\infty}|<C.
$$
\end{deplemma}

\begin{deplemma}
Assume that $\delta>0$ and there exist constants
$f_{\infty},\gamma_{1},\gamma_{2},\gamma_{3},C_{1}$,
$\tilde{C},C_{2},C_{3}>0,$ such that 
$$
\max_{|z|=n}|f_{n}(z)|\leq C_{1}(e^{\gamma_{1}\sqrt{n}}),
$$
$$
\max_{|z/n-1|<n^{-\alpha}}|f_{n}(z)|e^{-\frac{\gamma_{2}|z-n|}{\sqrt{n}}}\leq
C_{2},
$$
$$
\max_{|z/n-1|<n^{\delta-1}}|f_{n}(z)-f_{\infty}|e^{-\frac{\gamma_{3}|z-n|}{\sqrt
{n}}}\leq C_{3}.
$$
Let $a_{n}$ be a sequence of positive numbers, $|a_{n}|<\tilde{C}$,
and assume that
$$
\max_{|z/n-1|<n^{\delta-1}}|f_{n}(z)-f_{\infty}|e^{\frac{-\gamma_{1}|z-n|}{\sqrt
{n}}}\leq C_{1}a_{n}.
$$
Then there exists a constant
$C=C(\gamma_{1},\gamma_{2},\gamma_{3},C_{1},\tilde{C},C_{2},C_{3})$
such that for all ${n>0}$ we have
$$
|f_{nn}-f_{\infty}|<Ca_{n}.
$$
\end{deplemma}

\section{Asymptotics of the poissonized variance}
\label{sec_poiss}
\subsection{Estimates of the Bessel kernel}
 To prove Proposition \ref{poissmain} we need the following estimates of the
discrete Bessel kernel (we postpone the proof to the last Section).
Set 
$$
u_x=x/\sqrt{n},u_y=y/\sqrt{n},\phi_x=\arccos(u_x/2),\phi_y=\arccos(u_y/2).
$$
\begin{proposition}[points are far from the edge of the spectrum]
\label{bulk}
Let 
$$
|x|<2\sqrt{n}-n^{\delta_{1}},\,|y|<2\sqrt{n}-n^{\delta_{2}},\,\delta_{1},\delta_
{2}>1/6.
$$
Then
$$
|J(x,y,\theta^2)|\leq\frac{C\exp(\gamma|\theta-\sqrt{n}|)}{|e^{i\phi_{x}}-e^{
i\phi_{y}}|\sqrt[4]{2-u_{x}}\sqrt[4]{2-u_{y}}\sqrt{n}}+
$$
$$
+\exp\left(-n^{c(-1/6+\max(\delta_{1},\delta_{2}))}+\gamma|\theta-\sqrt{n}
|\right).
$$
\end{proposition}

\begin{proposition}[one of the points is at the edge of the spectrum]
\label{edge}
Let 
$$
|x|<2\sqrt{n}-n^{\delta},2\sqrt{n}-n^{1/6}<|y|<2\sqrt{n}+n^{1/6},\delta>1/6.
$$
Then
$$
|J(x,y,\theta^2)|\leq\left(\frac{C}{n^{5/12}\left(2-u_{x}\right)^{3/4}}+O(e^{-n^
{c(\delta-1/6)}})\right)\exp(\gamma|\theta-\sqrt{n}|).$$
\end{proposition}

\begin{proposition}[one of the points is beyond the edge of the spectrum]
\label{beyond}
Let
$$
|x|<2\sqrt{n},|y|>2\sqrt{n}+n^{\delta},\delta>1/6.
$$
Then
$$
|J(x,y,\theta^2)|\leq\left(Ce^{-n^{c(\delta-1/6)}}
\right)\exp(\gamma|\theta-\sqrt{n}|).
$$
\end{proposition}
\subsection{Proof of Proposition \ref{poissmain}}
We now rewrite the estimate from Proposition \ref{bulk}. 
$$
\frac{C}{|e^{i\phi_{x}}-e^{i\phi_{y}}|\sqrt[4]{2-u_{x}}\sqrt[4]{2-u_{y}}\sqrt{n}
}\leq\frac{C}{|\sqrt{2-u_{x}}-\sqrt{2-u_{y}}|\sqrt[4]{2-u_{x}}\sqrt[4]{2-u_{y}}
\sqrt{n}},
$$
Using the equality
$$
(a^{2}-b^{2})ab=(a^{4}-b^{4})\frac{ab}{a^{2}+b^{2}},
$$
we get
$$
J(x,y,\theta^2)\leq\frac{C\exp{\gamma|\theta-\sqrt{n}|}}{|x-y|}\left(\frac{\sqrt
[4]{2-u_{x}}}{\sqrt[4]{2-u_{y}}}+\frac{\sqrt[4]{2-u_{y}}}{\sqrt[4]{2-u_{x}}}
\right).
$$
Divide the interval of summation into three parts:
$$
R_{1}=\{(i,\, j)\in\mathbb{Z}^{2}|\, i\in(x+\delta,\, y-\delta),\, j\notin[x,\,
y]\}
$$
$$
R_{2}=\{(i,\, j)\in\mathbb{Z}^{2}|\,(i\in[x,\, x+\delta],\,
j\notin[i-\delta-1,\, x-1]\cup[x,y])\vee$$
$$
\vee(i\in[y-\delta,\, y],\, j\notin[y+1,\, i+\delta+1]\cup[x,y])\}
$$
$$
M=\{(i,\, j)\in\mathbb{Z}^{2}|\,(i\in[x,\, x+\delta],\, j\in[i-\delta-1,\,
x-1])\vee
$$
$$
\vee(i\in[y-\delta,\, y],\, j\in[y+1,\, i+\delta+1])\}
$$
where $\delta=\varepsilon(x-y)$ will be chosen later. Following \cite{key-3} we
are going to
estimate the sum from parts $R_{1}$ and $R_{2}$, which may be neglected, and
compute it in
\emph{M}. 

\paragraph{Estimate in $R_{1}$:}

Write down the sum from $R_1$:
\begin{multline}
\sum_{i\in(x+\delta,y-\delta)}\sum_{j\notin[x,y]}J^{2}(i,j,\theta^2)=\sum_{
i\in(x+\delta,y-\delta)}\sum_{j=-\infty}^{x-1}J^{2}(i,j,\theta^2)+\\
+\sum_{i\in(x+\delta,y-\delta)}\sum_{j=y+1}^{\infty}J^{2}(i,j,\theta^2).
\end{multline}
We estimate the first summand (the second one may be estimated similarly)

Summing in $j$ we have:
\begin{multline}
\label{R_1-1}
\sum_{j=y+1}^{\infty}J^{2}(i,j,\theta^2)=\sum_{j=y+1}^{2\sqrt{n}-n^{1/6}}J^{2}(i
,j,\theta^2)+\\
+\sum_{j=2\sqrt{n}-n^{1/6}}^{2\sqrt{n}+n^{1/6}}J^{2}(i,j,\theta^2)+\sum_{
j=2\sqrt{n}+n^{1/6}}^{\infty}J^{2}(i,j,\theta^2).
\end{multline}

We estimate the first summand from (\ref{R_1-1}):
\begin{multline}
\exp(-\gamma|\theta-\sqrt{n}|)\sum_{j=y+1}^{2\sqrt{n}-n^{1/6}}J^{2}(i,j,
\theta^2)\leq\\
\leq\sum_{j=y+1}^{2\sqrt{n}-n^{1/6}}\frac{C}{(j-i)^{2}}\left(\frac{\sqrt[4]{2-u_
{i}}}{\sqrt[4]{2-u_{j}}}+\frac{\sqrt[4]{2-u_{j}}}{\sqrt[4]{2-u_{i}}}\right)^{2}
\leq \\
\leq\sum_{j=y+1}^{2\sqrt{n}-n^{1/6}}\frac{C}{(j-i)^{2}}\left(\frac{\sqrt{2-u_{i}
}}{\sqrt{2-u_{j}}}+\frac{\sqrt{2-u_{j}}}{\sqrt{2-u_{i}}}\right)\leq \\
\end{multline}
\begin{multline}
\leq
C\intop_{y+1}^{2\sqrt{n}-n^{1/6}}\frac{1}{(t-i)^{2}}\left(\frac{\sqrt{2\sqrt{n}
-i}}{\sqrt{2\sqrt{n}-t}}+\frac{\sqrt{2\sqrt{n}-t}}{\sqrt{2\sqrt{n}-i}}
\right)dt=\\
=-C\frac{1}{t-i}\frac{\sqrt{2\sqrt{n}-t}}{\sqrt{2\sqrt{n}-i}}\bigg|_{y+1}^{
2\sqrt{n}-n^{1/6}}=\\
=-C\left(\frac{n^{1/6}}{(2\sqrt{n}-i-n^{1/6})\sqrt{2\sqrt{n}-i}}-\frac{\sqrt{
2\sqrt{n}-y-1}}{(y+1-i)\sqrt{2\sqrt{n}-i}}\right)\leq\\
\leq\frac{C}{y-i+1}.
\end{multline}
Now sum in $i$:
$$
\sum_{i\in(x+\delta,y-\delta)}\frac{C}{y-i+1}\leq C\ln(\frac{y-x}{\delta}).
$$
Now estimate the second summand:
$$
\sum_{j=2\sqrt{n}-n^{1/6}}^{2\sqrt{n}+n^{1/6}}J^{2}(i,j,\theta^2)\leq\sum_{
j=2\sqrt{n}-n^{1/6}}^{2\sqrt{n}+n^{1/6}}\frac{C\exp(\gamma|\theta-\sqrt{n}|)}{n^
{5/6}\left(2-u_{i}\right)^{3/2}}=C\frac{\exp(\poissdev)n^{1/12}}{(2\sqrt{n}-i)^{
3/2}}.
$$
Sum this estimate in $i$:
\begin{multline}
\sum_{i\in(x+\delta,y-\delta]}\frac{n^{1/12}}{(2\sqrt{n}-i)^{3/2}}\leq
Cn^{1/12}\intop_{x}^{y}\frac{1}{(2\sqrt{n}-t)^{3/2}}dt=\\
=Cn^{1/12}\left(\frac{1}{\sqrt{2\sqrt{n}-t}}\right)\bigg|_{x}^{y}\leq\frac{Cn^{
1/12}}{\sqrt{2\sqrt{n}-y}}\leq\\
\leq Cn^{1/12-1/4}=o(1).
\end{multline}
The third summand is exponentially small by Proposition \ref{beyond}.

\paragraph{Estimate in $R_{2}$:}
$$
\sum_{i\in[y-\delta,\, y]}\sum_{j\notin[y+1,\,
i+\delta+1]}J^{2}(x,y,\theta^2)\leq\sum_{i\in[y-\delta,\,
y]}\frac{\exp(\poissdev)}{\delta}=O(1)\exp(\poissdev).$$
\paragraph{Asymptotics in $M$:}
From the Debye asymptotics for the Bessel function we get:
\begin{equation}
\label{Debye}
J(x,y,\theta^2)=\frac{\sin(\phi_y(x-y))}{\pi(x-y)}+o(\exp(\poissdev)),
\end{equation}
(see the proof of Lemma 3.5 in \cite{key-2}) which implies
\begin{multline}
\sum_{i\in[y-\delta,\, y]}\sum_{j\in[y+1,\, i+\delta+1]}J^{2}(x,y,\theta^2)=\\
=\sum_{i\in[y-\delta,\, y]}\sum_{j\in[y+1,\,
i+\delta+1]}\frac{1}{\pi^{2}}\frac{\sin^{2}(\phi_y(i-j))}{(i-j)^{2}}
(1+o(1))+o(e^{\gamma|n-\theta^2|/\sqrt{n}}),
\end{multline}
In sums of this kind replacement of $\sin^2(\cdot)$ by its mean value
$\frac{1}{2}$ does not change the asymptotics. More precisely, we now use Lemma
4.6 from \cite{key-3} on our kernel making the change in parameters as follows: $$
i=\xi_{1}+y-\delta,\, j=\xi_{2}+y+1.
$$
We obtain
\begin{multline}
\sum_{i\in[y-\delta,\, y]}\sum_{j\in[y+1,\,
i+\delta+1]}\frac{1}{\pi^{2}}\frac{\sin^{2}(\phi(i-j))}{(i-j)^{2}}=\\
=\sum_{i\in[y-\delta,\, y]}\sum_{j\in[y+1,\,
i+\delta+1]}\frac{1}{2\pi^{2}}\frac{1}{(i-j)^{2}}(1+o(1))=\\
=\frac{1}{2\pi^{2}}\ln\delta(1+o(1)).
\end{multline}
Now set $\delta=\dfrac{|x-y|}{\ln(|x-y|)}$. This choice of $\delta$ implies that
the sum in $R_1\leq C\ln\ln(x-y)$, and that finishes the proof of Proposition
\ref{poissmain}.

\section{Depoissonization}
\label{sec_depoiss}
In this section we finish the proof of Theorem \ref{main}. Note that the
depoissonization techniques described above may be applied only to quantities
linear in measure (i.e. to an expectation, not variance). That is why we will now,
following \cite{key-1},  find an expectation with the same asymptotics as our
poissonized and original variances has.
To do this we will need Lemma 6.3 from \cite{key-1} and its depoissonization:
\begin{proposition}
\label{dep-1}
There exists $\varepsilon_0>0$ such that the following holds. For all
$\delta_0>\dfrac{1}{6}$ there exist constants $C>0$, $\gamma>0$ such that for
all $n\in\mathbb{N}$, all $x,|x|<2\sqrt{n}-n^{\delta_0}$ and all
$\theta\in\mathbb{C}$ such that
$$
\left|\frac{\theta}{\sqrt{n}}-1\right|<\varepsilon_0
$$
we have
$$
|J(x,x,\theta^2)-\frac{1}{\pi}\arccos\frac{x}{2\sqrt{n}}|\leq\frac{
\exp(\gamma|\theta-\sqrt{n}|)}{2\sqrt{n}-x}.
$$
\end{proposition}
\begin{proposition}
\label{dep-2}
For all $\delta_0>\frac{1}{6}$ there exists a constant $C$ such that
$$
|\mathbb{E}_{\mathbb{P}l^{(n)}}(c_x)-\frac{1}{\pi}\arccos\frac{x}{2\sqrt{n}}
|\leq\frac{C(2-|u_{x}|)}{\sqrt{n}},
$$
for all $x,|x|<2\sqrt{n}-n^{\delta_0}$.
\end{proposition}
Write
\begin{multline}
Var_{\theta^2}\left(\sum_{i=k}^{k+l}c_{i}\right)=\mathbb{E}_{\theta^2}
\left(\sum_{i=k}^{k+l}c_{i}-\sum_{i=k}^{k+l}\mathbb{E}_{\theta^2}c_{i}\right)^{2
}=\\
=\mathbb{E}_{\theta^2}\left(\sum_{i=k}^{k+l}c_{i}-\sum_{i=k}^{k+l}J(i,i,
\theta^2)\right)^{2}.
\end{multline}
Set
$$
\Omega(t)=\begin{cases}
\frac{2}{\pi}(t\arcsin(t/2)+\sqrt{4-t^{2}}), & \text{if
\ensuremath{|t|\leq2};}\\
|t|, & \text{if \ensuremath{|t|>2}, }\end{cases}
$$
$$
F_{\lambda}(t)=\Phi_{\lambda}(t)-\sqrt{n}\Omega(t/\sqrt{n}),
$$
where $\Phi_{\lambda}(t)$ is a function representing the upper edge of the
diagram $\lambda$. Note that $\Omega(t)$ is the limit shape of Plancherel Young
diagrams (see \cite{key-7}, \cite{key-6} for details).
Now write
\begin{multline*}
F_{\lambda}(k+1)-F_{\lambda}(k)=1-2c_{k}(\lambda)-\sqrt{n}
\left(\Omega\left(\frac{k+1}{\sqrt{n}}\right)-\Omega\left(\frac{k}{\sqrt{n}}
\right)\right)=\\
=2\left(\frac{\arccos\frac{k}{2\sqrt{n}}}{\pi}-c_{k}(\lambda)\right)+\frac{2}{
\pi}\arcsin\left(\frac{k}{2\sqrt{n}}\right)-\sqrt{n}\left(\Omega\left(\frac{k+1}
{\sqrt{n}}\right)-\Omega\left(\frac{k}{\sqrt{n}}\right)\right)=\\
=2\left(J(k,k;\,\theta^{2})-c_{k}(\lambda)\right)+2\left(\frac{\arccos\frac{k}{
2\sqrt{n}}}{\pi}-J(k,k;\,\theta^{2})\right)+\\
+\left(\,\frac{2}{\pi}\arcsin\left(\frac{k}{2\sqrt{n}}\right)-\sqrt{n}
\left(\Omega\left(\frac{k+1}{\sqrt{n}}\right)-\Omega\left(\frac{k}{\sqrt{n}}
\right)\right)\right).\end{multline*}
From the Taylor formula for $\Omega$ we have
$$
\left|\sqrt{n}\left(\Omega\left(\frac{k+1}{\sqrt{n}}\right)-\Omega\left(\frac{k}
{\sqrt{n}}\right)\right)-\frac{2}{\pi}\arcsin\left(\frac{k}{2\sqrt{n}}
\right)\right|\le\frac{10}{\sqrt{4n-k^{2}}}.
$$
Summing these estimates and applying Proposition \ref{dep-1} we obtain that
after the depoissonization
$$
\lim_{n\rightarrow\infty}\frac{\mathbb{E}_{\mathbb{P}l^{(n)}}\left(F_{\lambda}
(k+l)-F_{\lambda}(k)\right)^{2}}{\log l}=\frac{1}{\pi^{2}}.
$$
In the same way, applying Proposition \ref{dep-2}, we finish the proof of
Theorem \ref{main}.
\section{Proofs of the estimates of the Bessel kernel}
\label{sec_est}
\paragraph{Plan of the proofs}
We are using the Okounkov's integral representation of the kernel (\cite{key-4}):
$$
J(x,y,\theta^{2})=\frac{1}{(2\pi
i)^{2}}\intop\intop_{|z|<|w|}\frac{\exp(\theta\left(z-\frac{1}{z}-w+\frac{1}{w}
\right))}{(z-w)z^{x+1}w^{-y}}dzdw.
$$

Set
$$
u_{x}=\frac{x}{\sqrt{n}},\, S(z,u)=z-\frac{1}{z}-u\ln
z,\,\phi_{x}=\arccos\frac{u_{x}}{2}.
$$
Denote
$$
\Phi(z,w,\theta)=\frac{\exp\theta\left(z-\frac{1}{z}-w+\frac{1}{w}\right)}{
(z-w)z^{x+1}w^{-y}}.
$$
Now note that
$$
|\Phi(z,w,\theta)|\leq\exp(\poissdev)|\Phi(z,w,\sqrt{n}
)|=\exp(\poissdev)\left|\frac{\exp\sqrt{n}(S(z,u_x)-S(z,u_y))}{z(z-w)}\right|.
$$
This allows us to work with the Bessel kernel of a real parameter,

The proof is carried out using the steepest descent method. Following \cite{key-4}
we deform initial contours so that they pass through critical points of $S(z,u)$
and then we estimate the absolute value of the integral splitting contours in
parts.

\subsection{Proof of Proposition \ref{bulk}}
Without a loss of generality assume that $x>y.$ Set
$I_{x}=[-n^{-\beta_{x}},n^{-\beta_{x}}],I_{y}=[-n^{-\beta_{y}},n^{-\beta_{y}}]$,
$\beta_{x}$ and $\beta_{y}$ will be defined later. 

\begin{figure}[h]
\includegraphics{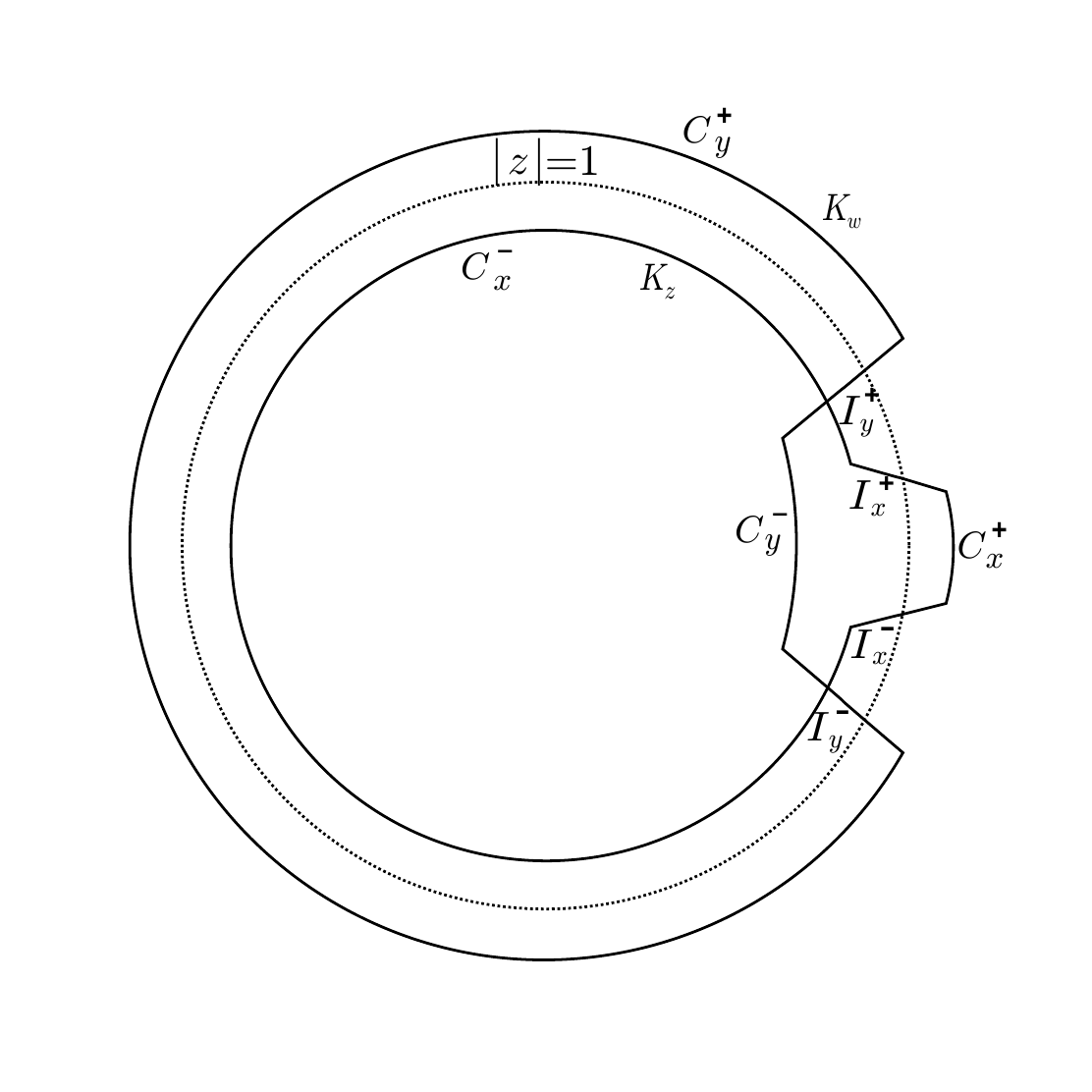}
\caption{deformed contours.}
\label{fig_contour}
\end{figure}

We now deform the contours in the following way. Each one consists of four
parts: two circular arcs and two intervals transversal to the unit circle,
similarly to \cite{key-1} (see Fig. \ref{fig_contour}). 

Introduce the parametrization on the intervals: 
$$
I_{x}^{+}(t)=e^{i\phi_{x}}+At,\, I_{x}^{-}(t)=e^{-i\phi_{x}}+At,
$$
$$
I_{y}^{+}(t)=e^{i\phi_{y}}+Bs,\, I_{y}^{-}(t)=e^{-i\phi_{y}}+Bs,$$
$$
t\in I_{x},s\in I_{y}.
$$
Here $A,B\in\mathbb{C}$ are independent on $n$ and chosen such that the
intervals are transversal to the unit circle.

Let $C_{x}^{\pm},C_{y}^{\pm}$ be the circular arcs of the contours (- for the inner
part and + for the outer),
0 be the center of both circles and $1\pm cn^{-\beta_{x}}$, $1\pm
cn^{-\beta_{y}}$ be their radii.

Write
$$
S(z(t),u)=e^{i\phi_{x}}+At-\frac{1}{e^{i\phi_{x}}+At}-u_x\ln\left(e^{i\phi_{x}}
+At\right),
$$
on the interval.
$$
\frac{dS}{dt}\bigg|_{t=0}=\left(A+\frac{A}{\left(e^{i\phi_{x}}+At\right)^{2}}
-\frac{u_xA}{e^{i\phi_{x}}+At}\right)\bigg|_{t=0}=0,
$$
\begin{multline}
\frac{d^{2}S}{dt^{2}}\bigg|_{t=0}=\left(-\frac{2A^{2}}{\left(e^{i\phi_{x}}
+At\right)^{3}}+\frac{u_xA^{2}}{\left(e^{i\phi_{x}}+At\right)^{2}}\right)\bigg|_
{t=0}=A^{2}(u_xe^{i\phi_{x}}-2)e^{-3i\phi_{x}}=\\
=A^{2}e^{-3i\phi_x}(u_x\cos\phi_{x}+iu_x\sin\phi_{x}-2)=2A^{2}e^{-3i\phi_{x}}
\sin\phi_{x}(i\cos\phi_{x}-\sin\phi_{x})=\\
=2iA^{2}e^{-2i\phi_{x}}\sin\phi_{x}.
\end{multline}
We used the fact that $u_x=2\cos\phi_x$. We obtain
$$
\left|\Re\frac{d^{2}S}{dt^{2}}\bigg|_{t=0}\right|\geq C\sqrt{2-u}.
$$
Recall that $f=\Theta(g) \iff f=O(g) \text{ and } g=O(f)$.
\begin{multline}
\frac{d^{3}S}{dt^{3}}\bigg|_{t=0}=\left(\frac{6A^{3}}{\left(e^{i\phi_{x}}
+At\right)^{4}}-\frac{2u_xA^{3}}{\left(e^{i\phi_{x}}+At\right)^{3}}\right)|_{t=0
}=6A^{3}e^{-4i\phi_{x}}-2u_xA^{3}e^{-3i\phi_{x}}=\\
=2A^{3}e^{-4i\phi_{x}}(3-u_xe^{i\phi_{x}})=\Theta(1).
\end{multline}
\begin{multline}
\frac{d^{4}S}{dt^{4}}\bigg|_{t=0}=\left(-\frac{24A^{4}}{\left(e^{i\phi_{x}}
+At\right)^{5}}+\frac{6u_xA^{4}}{\left(e^{i\phi_{x}}+At\right)^{4}}\right)|_{t=0
}=24A^{4}e^{-5i\phi_{x}}-6u_xA^{4}e^{-4i\phi_{x}}=\\
=6A^{4}e^{-5i\phi_{x}}(-4+u_xe^{i\phi_{x}})=\Theta(1).
\end{multline}
Derivatives $I_{y}^{\pm}(s)$ are estimated similarly.

\paragraph{Circular parts.}
\label{circ_parts}
Switch to polar coordinates $z=(1+t)e^{i\phi}$.
\begin{multline}
\Re S(z,u)=(1+t-1+t-t^{2})\cos\phi(t)-u(t-t^{2}/2)+O(t^{3})=\\
=(2\cos\phi(t)-u)(t-t^{2}/2)+O(t^{3})=(2\cos(\arccos\frac{u}{2}+c_{1}t)-u)(t-t^{
2}/2)+O(t^{3})=\\
=2ct(-\sin\arccos\frac{u}{2})t+O(t^{3})=-2c\sqrt{1-\frac{u^{2}}{4}}t^{2}+O(t^{3}
)=\\
=-c\sqrt{2-u}\sqrt{2+u}t^{2}+O(t^{3}),
\end{multline}
$$
\sqrt{n}\Re S(z,u)\leq-c\sqrt{n}\sqrt{2-u}t^{2}+\sqrt{n}O(t^{3}).
$$

We want to deform the contours so that their intersection points are far enough
from the unit circle, where $\Phi(z,w,\theta)$ is exponentially small. To
guarantee that, we need the following conditions. Set $t=n^{-\beta}$ on both ends
of the interval and write

1)  $\sqrt{n(2-u)}t^{2}\to\infty$ as $n\to\infty$,

2) $\dfrac{t^{3}}{\sqrt{2-u}t^{2}}\to0$ as $n\to\infty$,

3) $n^{-\beta}\leq c|e^{i\phi_{x}}-e^{i\phi_{y}}|$

Now we verify these conditions. $\sqrt{2-u}\geq n^{\delta/2-1/4},$ so 1) gives
$$
\delta/2+1/4-2\beta>0\implies\beta<\delta/4+1/8,
$$
and 2) gives
$$
\delta/2-1/4>-\beta\implies\beta>1/4-\delta/2.
$$
$\beta$ we need exists if
$$
\delta/4+1/8>1/4-\delta/2\implies3\delta/4>1/8,\,\delta>\frac{1}{6}.
$$
Now consider 3). Let 
$$
x=2\sqrt{n}-n^{\delta}+an^{\alpha},\,
y=2\sqrt{n}-n^{\delta}+bn^{\alpha},\,\delta>1/6,\,\alpha<\delta.
$$

\begin{figure}[h]
\includegraphics[scale=0.8]{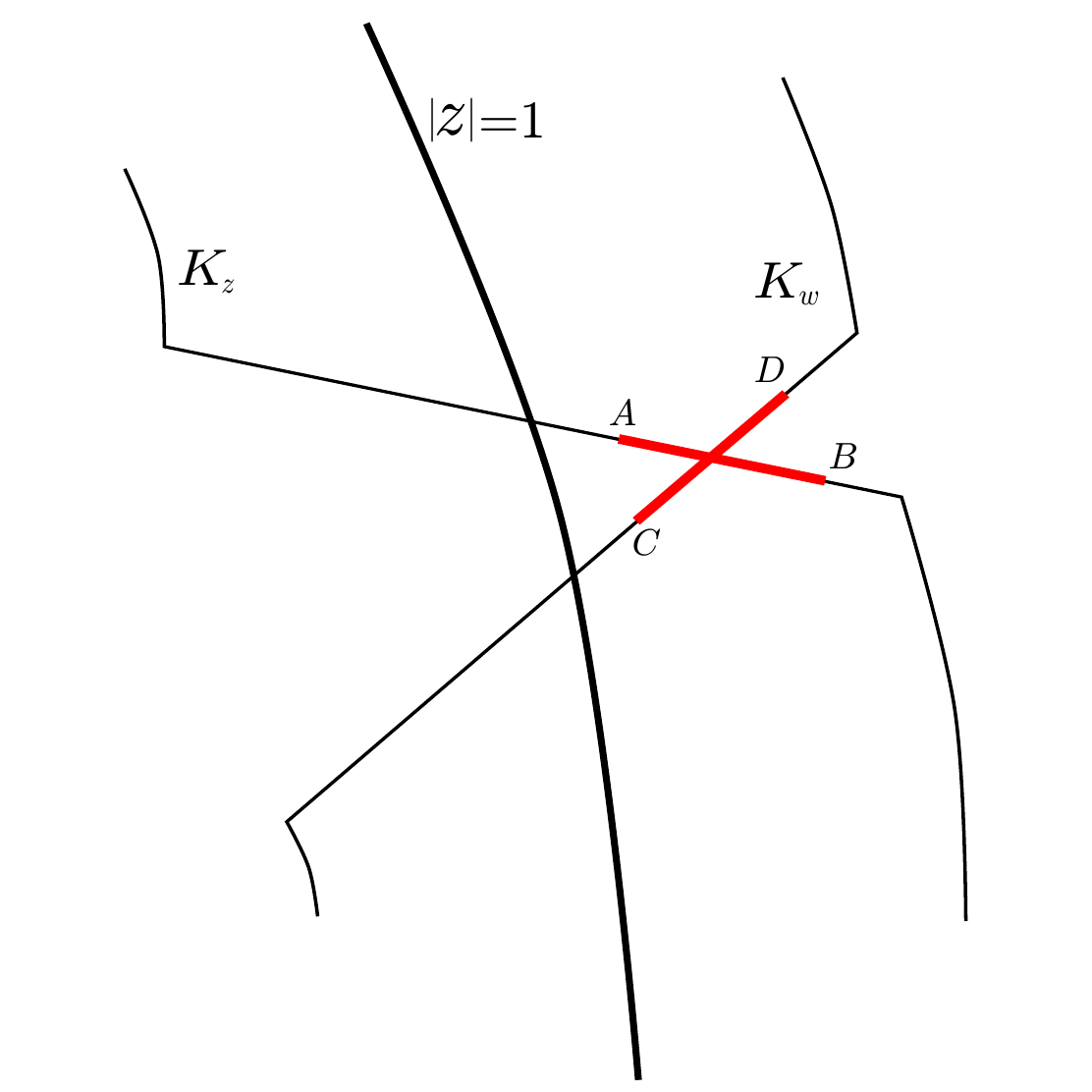}
\caption{The integral along $AB$ and $CD$ gives the main contribution to the
estimate.}
\label{fig_bulk}
\end{figure}

Then
$$
|e^{i\phi_{x}}-e^{i\phi_{y}}|=\Theta(|\sin\phi_{x}-\sin\phi_{y}|)=\Theta(|\sqrt{
2-u_{x}}-\sqrt{2-u_{y}}|),
$$
$$
\sqrt{n^{\delta-1/2}+an^{\alpha-1/2}}-\sqrt{n^{\delta-1/2}+bn^{\alpha-1/2}}=n^{
\delta/2-1/4}\left(\sqrt{1+an^{\alpha-\delta}}-\sqrt{1+bn^{\alpha-\delta}}
\right)=
$$
$$
=\Theta(n^{\alpha-1/4-\delta/2}).
$$
as $n\to\infty$, and in this case there exists $C>0$ such that the following
inequality holds:
$$
\frac{C}{|e^{i\phi_{x}}-e^{i\phi_{y}}|\sqrt[4]{2-u_{x}}\sqrt[4]{2-u_{y}}\sqrt{n}
}\geq\frac{C}{(a-b)n^{\alpha}}.
$$
Rewrite 3) as follows:
$$
n^{-\delta/4-1/8}<n^{-\beta}<n^{\alpha-1/4-\delta/2},
$$
$$
-\delta/4-1/8<\alpha-1/4-\delta/2,
$$
$$
\alpha>1/8+\delta/4.
$$
It remains to consider
$$
\alpha<1/8+\delta/4.
$$
Note that in this case 
$$
\alpha-1/4-\delta/2<-\alpha.
$$
We now estimate the integral in a neighborhood of the intersection along the
intervals $J_{1},J_{2}$ of length 
$$
n^{\alpha-1/4-\delta/2}
$$
(the intervals $AB$ and $CD$ on Fig. \ref{fig_bulk}).

\begin{multline}
\left|\intop_{J_{1}}\intop_{J_{2}}\Phi(z,w,\theta)dzdw\right|\leq\intop_{J_{1}}
\intop_{J_{2}}|\Phi(z,w,\theta)|dzdw\leq\\
\leq\exp(\poissdev)\intop_{J_{1}}\intop_{J_{2}}|\Phi(z,w,\sqrt{n}
)|dzdw\leq\exp(\poissdev)\intop_{J_{1}}\intop_{J_{2}}\left|\frac{1}{\sqrt{c_{1}
s^{2}+c_{2}t^{2}}}\right|dsdt\leq\\
\leq \exp(\poissdev)C\max(|J_{1}|,|J_{2}|)\leq Cn^{-\alpha}.
\end{multline}

\paragraph{Main contribution.}
\begin{multline}
\label{int1}
\left|\intop_{I_{x}^{+}}\intop_{I_{y}^{+}}\Phi(x,y,
\theta)dzdw\right|\leq\exp(\poissdev)\intop_{I_{x}}\intop_{I_{y}}\left|\frac{
\exp(\sqrt{n}\left(S\left(z(t),u_{x}\right)-S(w(t),u_{y}\right)))}{e^{i\phi_{x}}
-e^{i\phi_{y}}+At-Bs}\right|dsdt\leq\\
\leq\exp(\poissdev)\frac{\intop_{I_{x}}\intop_{I_{y}}|\exp(\sqrt{n}
\left(S\left(z(t),u_{x}\right)-S\left(w(t),u_{y}\right)\right))|dsdt}{|e^{i\phi_
{x}}-e^{i\phi_{y}}|}=\\
=\exp(\poissdev)\frac{\intop_{I_{x}}\intop_{I_{y}}|\exp\sqrt{n}\left(S_{x}
''(0)t^{2}-S_{y}''(0)s^{2}+S_{x}^{(3)}(0)t^{3}-S_{y}{}^{(3)}(0)s^{3}+o(s^{3},t^{
3})\right)|dsdt}{|e^{i\phi_{x}}-e^{i\phi_{y}}|}.
\end{multline}
Write
$$
s=\frac{s'}{\sqrt[4]{n}\sqrt{|S_{y}''(0)|}},t=\frac{t'}{\sqrt[4]{n}\sqrt{|S_{x}
''(0)|}}.
$$

Note that one can choose $\beta$ so that after change of variables the coefficient in
the third term of the Taylor expansion goes to 0 as $n \to \infty$. Indeed, length
of the new intervals of integration is equal to 
$$
n^{-\beta+1/4}\sqrt[4]{2-u}
$$
and
\begin{multline}
\frac{\sqrt{n}t^{'3}}{n^{3/4}\left(\sqrt{|S_{x}''|}\right)^{3}}\leq
Cn^{-1/4}n^{3/8-3\delta/4}n^{-3\beta+3/4+3\delta/4-3/8}=Cn^{1/2-3\beta}\leq\\
\leq Cn^{1/2-3(1/8+\delta/4-\varepsilon)}=Cn^{1/8-3\delta/4+3\varepsilon}.\\
\end{multline}
We get
\begin{multline}
\left|\intop_{I_{x}^{+}}\intop_{I_{y}^{+}}\Phi(x,y,\theta)dzdw\right|\leq
C\exp(\poissdev)\frac{|e^{i\phi_{x}}-e^{i\phi_{y}}|^{-1}}{\sqrt{n}\sqrt{S_{x}''}
\sqrt{S_{y}''}}\leq\\
\leq
C\exp(\poissdev)\frac{|e^{i\phi_{x}}-e^{i\phi_{y}}|^{-1}}{\sqrt{n}\sqrt[4]{2-u_{
x}}\sqrt[4]{2-u_{y}}}.
\end{multline}

Same computations give the estimates of
$|\intop_{I_{x}^{-}}\intop_{I_{y}^{-}}\Phi(x,y,\theta)dzdw|,\,|\intop_{I_{x}^{+}
}\intop_{I_{y}^{-}}\Phi(x,y,\theta)dzdw|,$
$|\intop_{I_{x}^{-}}\intop_{I_{y}^{+}}\Phi(x,y,\theta)dzdw|$ and the estimate
of the
integral along the linear intervals in the case $\alpha\leq1/8+\delta/4$. 

\paragraph{Residue.}

As we deform the contours we pick up the residue at $z=w$:
$$
|Res|=\left|\dfrac{1}{2\pi
i}\int_{re^{-i\phi}}^{re^{i\phi}}\dfrac{1}{w^{x-y+1}}dw\right|=\left|\dfrac{1}{
2\pi ir^{x-y}}{\displaystyle
\int_{e^{-i\phi}}^{e^{i\phi}}\dfrac{1}{t^{x-y+1}}dt}\right|=\left|\dfrac{1}{2\pi
ir^{x-y}}K_{\sin}(x-y,\phi)\right|.
$$
Note that $sgn(x-y)=sgn(r-1).$
\subsection{Proof of Proposition \ref{edge}}
\begin{figure}[h]
\includegraphics[scale=0.8]{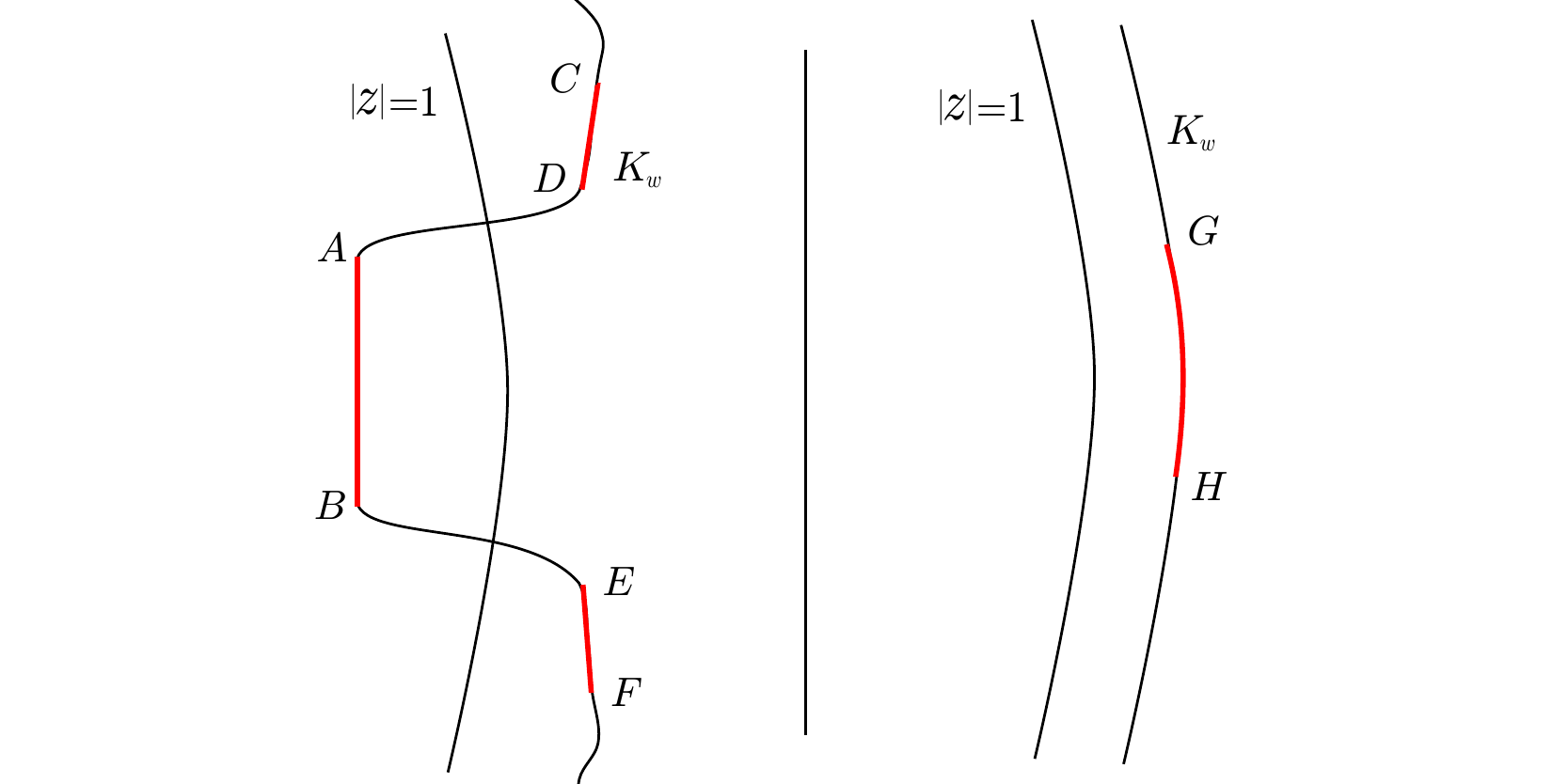}
\caption{deformation of the contours for $y<2\sqrt{n}$ and $y>2\sqrt{n}$
respectively.}
\label{fig_edge}
\end{figure}

We deform the contours similarly to the case of Proposition
\ref{bulk}, with $\beta_{y}=1/6$. In the case $y>2\sqrt{n}$ let contour $K_{w}$
be a circle with center in 0 and of radius $1+n^{-\beta_{y}}$. Consider parts of
contours of length at most 
$$
\sqrt{|2-u_{y}|}<n^{-1/6}.
$$ adjacent to  $I_x^\pm$ and $I_y^\pm$ (for $y>2\sqrt{n}$ consider a circular part
in the neighborhood of real axis) ($AB,CD,EF$ and $FG$ of Fig. \ref{fig_edge}).

Write
$$
\Re S(z,u)=(2\cos\phi(t)-u)(t-t^{2}/2)+O(t^{3}).
$$
On these parts $\exp(n^{1/2}\Re S(z,u))=O(1)$ and beyond them 
on circular parts $n^{1/2}\Re S(z,u)$ is large:
while $|2\cos\phi(t)-u|>\sqrt{|2-u_{y}|}$ we have
$$
n^{1/2}\Re S(z,u)<-n^{1/6}.$$
Now estimate the contribution from integration along the intervals $I_x^{\pm}$ and
$I_y^{\pm}$ similarly to (\ref{int1}) making another change of variables: 
$$
s=\frac{s'}{n^{1/6}},t=\frac{t'}{\sqrt[4]{n}\sqrt{|S_{x}''|}},$$
and
$$
|z-w|>C|\sqrt{2-u_{x}}|.
$$
The residue and contribution from the circular parts are estimated similarly to the
proof of Proposition \ref{bulk}.

\subsection{Proof of Proposition \ref{beyond}}
In this case one may not deform the contours choosing their radii so that
$\Phi(z,w,\theta)$ is exponentially small, see the estimate of the contribution
from the circular parts in Proposition \ref{bulk}.

\section{General local patterns}
\label{corners}
\subsection{Variance for general local patterns}
Denote $\vec x=\{x_1,...,x_k\} \subset \mathbb{Z}$. Denote $c_{i+\vec
x}=c_{i+x_1}c_{i+x_2}\cdots c_{i+x_k}$.
\begin{proposition}
 Let $\vec x \subset \mathbb{Z}$ be a finite set, $a,b \in (-2,2)$, and
${x_n},{y_n}$ be an integer sequence such that 
$$
\lim\limits_{n\to\infty}x_n/\sqrt{n}=a,\lim\limits_{n\to\infty}y_n/\sqrt{n}=b,
y_n>x_n
$$
and $\lim\limits_{n\to\infty}(y_n-x_n)=+\infty$. Then there exists a constant
$C>0$ such that
$$
Var_{\Pl^{(n)}}\left(\sum_{i=x_n}^{y_n}c_{i+\vec x}\right) \leq C(y_n-x_n).
$$
\end{proposition}
We start with the proof of the poissonized version of this proposition:
\begin{equation}
\label{poiss_conf}
Var_{\theta^2}\left(\sum_{i=x_n}^{y_n}c_{i+\vec x}\right) \leq
(y_n-x_n)C\exp(\gamma|\theta-\sqrt{n}|).
\end{equation}
Rewrite our variance in the following way:
\begin{multline}
\label{conf_var}
Var_{\theta^2}\left(\sum_{i=x_n}^{y_n}c_{i+\vec
x}\right)=\mathbb{E}_{\theta^2}\left(\sum_{i}c_{i+\vec
x}\right)+\mathbb{E}_{\theta^2}\left(\sum_{i \neq j}c_{i+ \vec x}c_{j+ \vec
x}\right)-\\
-\sum_{i}\left(\mathbb{E}_{\theta^2}c_{i+\vec x}\right)^2-\sum_{i \neq
j}\mathbb{E}_{\theta^2}c_{i+\vec x}\mathbb{E}_{\theta^2}c_{j+\vec x}.
\end{multline}
Estimate the sum on the right hand side of (\ref{conf_var}):
$$
\mathbb{E}_{\theta^2}\left(\sum_{i}c_{i+\vec
x}\right)-\sum_{i}\left(\mathbb{E}_{\theta^2}c_{i+\vec x}\right)^2 \leq
C_1(y_n-x_n)\exp(\poissdev);
$$
$$
\left|\sum_{i,j:i+\vec x \cap j+\vec
x\neq\varnothing}\mathbb{E}_{\theta^2}\left(\sum_{i \neq j}c_{i+ \vec x}c_{j+
\vec x}\right)-\sum_{i,j:i+\vec x \cap j+\vec
x\neq\varnothing}\mathbb{E}_{\theta^2}c_{i+\vec x}\mathbb{E}_{\theta^2}c_{j+\vec
x}\right| \leq C_2(y_n-x_n)\exp(\poissdev). 
$$
Note that if $i+\vec x \cap j+\vec x=\varnothing$, then
$$
\mathbb{E}_{\theta^2}\left(c_{i+ \vec x}c_{j+ \vec
x}\right)-\mathbb{E}_{\theta^2}c_{i+\vec x}\mathbb{E}_{\theta^2}c_{j+\vec x}\leq
\frac{C_3\exp(\poissdev)}{(i-j)^2},
$$
from the determinantal form of expectations.
Summing this inequality in $i$ and $j$ we get the poissonized proposition.

The depoissonization is carried out similarly to the depoissonization of Theorem
\ref{main}. That is (\ref{poiss_conf}), the Debye asymptotics of the Bessel kernel
\ref{Debye} and Proposition \ref{dep-1} yield

\begin{equation}
\label{moments}
\mathbb{E}_{\theta^2}\left(\sum_{i=x_n}^{y_n}c_{i+\vec
x}-\mathbb{E}_n\left(\sum_{i=x_n}^{y_n}c_{i+\vec x}\right)\right)^2 \leq
(y_n-x_n)C\exp{\gamma|\theta-\sqrt{n}|}.
\end{equation}

Note that in the equation (\ref{moments}) expression
$\mathbb{E}_n\left(\sum_{i=x_n}^{y_n}c_{i+\vec x}\right)$ is just a constant
that does not depend on $\theta$. After depoissonization we get
$$
\mathbb{E}_{\Pl^{(n)}}\left(\sum_{i=x_n}^{y_n}c_{i+\vec
x}-\mathbb{E}_n\left(\sum_{i=x_n}^{y_n}c_{i+\vec x}\right)\right)^2 \leq
C(y_n-x_n).
$$
Now using the Debye asymptotics we replace
$\mathbb{E}_n\left(\sum_{i=x_n}^{y_n}c_{i+\vec x}\right)$ by the expectation
with respect to the Plancherel measure and obtain the Proposition.
\paragraph{Lower bound}
We now give an example of a local pattern with a linearly growing variance.

\begin{proposition}
Assume that $\lim\limits_{n\to\infty}x_n/\sqrt{n} = 1,(y_n-x_n)=o(n^{1/6})>0$.
There
exist constants $c>0,n_0\in\mathbb{N}$ such that
$$
Var_{\Pl^{(n)}}\left(\sum_{i=x_n}^{y_n}c_ic_{i+1}\right)\geq c(y_n-x_n) 
$$
for $n>n_0$.
\end{proposition} 

Condition $(y_n-x_n)=o(n^{1/6})$ yields that we may compute our variance using the
sine kernel with a fixed parameter.
\begin{multline}
Var_{\Pl^{(n)}}\left(\sum_{i=x_n}^{y_n}c_ic_{i+1}\right)\sim\sum_i\mathbb{E}_{
\Pl^{(n)}}(c_ic_{i+1})-\sum_i\mathbb{E}_{\Pl^{(n)}}(c_ic_{i+1})^2+\\
+2\sum_i\mathbb{E}_{\Pl^{(n)}}(c_ic_{i+1}c_{i+2})-2\sum_i\mathbb{E}_{\Pl^{(n)}}
(c_ic_{i+1})\mathbb{E}_{\Pl^{(n)}}(c_{i+1}c_{i+2})+\\
+\sum_{|i-j|>1}\mathbb{E}_{\Pl^{(n)}}(c_ic_{i+1}c_jc_{j+1})-\sum_{|i-j|>1}
\mathbb{E}_{\Pl^{(n)}}(c_ic_{i+1})\mathbb{E}_{\Pl^{(n)}}(c_jc_{j+1}).
\end{multline}
Now note the following:
$K_{\sin}(0,\pi/2)=1/2,K_{\sin}(x,\pi/2)=0$ for $x=2k,k\neq0,$ and
$K_{\sin}(x,\pi/2)=\dfrac{(-1)^{(x-1)/2}}{x}$ for $x=2k+1$.
We get
\begin{multline}
 Var_{\Pl^{(n)}}\left(\sum_{i=x_n}^{y_n}c_ic_{i+1}\right)\sim\frac{1}{\pi^2}
\sum_i(\frac{\pi^2}{4}-1)-\frac{1}{\pi^4}\sum_i(\frac{\pi^2}{4}-1)^2+\\
+2\frac{1}{\pi^3}\sum_i(\frac{\pi^3}{8}-\pi)-2\frac{1}
{\pi^4}\sum_i(\frac{\pi^2}{4}-1)^2+\\
+\frac{1}{\pi^4}\sum_{|i-j|>1}\left(\frac{\sin^4(\frac{\pi}{2}(i-j))}{(i-j)^4}
+\frac{\sin^2(\frac{\pi}{2}(i-j+1))\sin^2(\frac{\pi}{2}(i-j-1))}{
(i-j-1)^2(i-j+1)^2}\right)-\\
-\frac{1}{\pi^4}\sum_{|i-j|>1}\left(2\frac{\sin^2(\frac{\pi}{2}(i-j))}{(i-j)^2}
(1+\pi^2/4)+2\frac{\sin(\frac{\pi}{2}(i-j+1))\sin(\frac{\pi}{2}(i-j-1))}{
(i-j-1)(i-j+1)}\right)-\\
-\frac{1}{\pi^4}\sum_{|i-j|>1}\frac{\pi^2}{4}\left(\frac{\sin^2(\frac{\pi}{2}
(i-j+1))}{(i-j+1)^2}+\frac{\sin^2(\frac{\pi}{2}(i-j-1))}{(i-j-1)^2}
\right)\sim\left(\frac{1}{12}-\frac{3}{8\pi^2}\right)(y_n-x_n).
\end{multline}
Here $f \sim g$ means $\lim\limits_{n\to\infty}\frac{f}{g}=1$.

Another interesting example of a local configuration is a corner: $c_i - c_i
c_{i+1}$. Note that from the formula
$$
Var(X+Y)=Var(X)+Var(Y)+Cov(X,Y)
$$
and the Cauchy-Bunyakovsky inequality we get a linear growth of the variance of the
corners' linear statistic in the same regime:

\begin{proposition}
\label{lv}
Assume that $\lim\limits_{n\to\infty}x_n = 1,(y_n-x_n)=o(n^{1/6})>0$. There
exist constants $c>0,n_0\in\mathbb{N}$ such that
$$
Var_{\Pl^{(n)}}\left(\sum_{i=x_n}^{y_n}(c_i-c_ic_{i+1})\right)\geq c(y_n-x_n) 
$$
for $n>n_0$.
\end{proposition}
 
\subsection{The central limit theorem} We now show that proposition \ref{lv}
yields the central limit theorem for the linear statistic of corners. Consider
the space of sequences $c_i \in \{0,1\}^{\mathbb{Z}}$ with the probabilistic measure
defined by the sine kernel. Let
$\vec{x}\subset\mathbb{Z}$ be a finite set. A process $c'_i = c_{i+\vec{x}}$ is a
stationary process with $Cov(c_i',c_{i+n}')$ decrease rate
$\Theta(\frac{1}{n^2})$. We now recall the conditions of the central limit
theorem from \cite{key-14}.

Let $X$ be a complete separable metric space, $\mathcal{F}$ a $\sigma$-algebra,
$P$ a probability measure, $T:X \to X$ a measurable map. Let also $T$ be
measure-preserving and dynamical system $(X,T,P)$ be ergodic. Define operator
$$
\hat{T}:L^2(X) \to L^2(X), \hat{T}\phi = \phi \circ T.
$$
and let $\hat{T}^*$ be its dual. Let $f$ be a random variable with
$\mathbb{E}(f)=0$, $\mathcal{F}_0$ be a sub-$\sigma$-algebra of $\mathcal{F}$
and set $\mathcal{F}_i=T^{-i}\mathcal{F}_0$.
\begin{theoremnn}[\cite{key-14}]
  If $\mathcal{F}_i$ is coarser than $\mathcal{F}_{i-1}$ and, for each $\phi \in
L^{\infty}(X)$, we have
$$
\mathbb{E}(\hat{T}\hat{T}^*\phi|\mathcal{F}_1)=\mathbb{E}(\phi|\mathcal{F}_1),
$$
then, for each $f\in L^{\infty}(X),\mathbb{E}(f)=0$ and
$\mathbb{E}(f|\mathcal{F}_0)=f$, such that
\begin{enumerate}
 \item $\displaystyle{\sum_{n=0}^{\infty}|\mathbb{E}(f\hat{T}^nf)|<\infty}$,
 \item the series
$\displaystyle{\sum_{n=0}^{\infty}\mathbb{E}(\hat{T}^{*n}f|\mathcal{F}_0)}$
converges absolutely almost surely,
\end{enumerate}
the sequence
$$
\frac{1}{\sqrt{n}}\sum_{i=0}^{n-1}\hat{T}^i f
$$
converges in distribution to a Gaussian random variable of zero mean and finite
variance.
\end{theoremnn}

To apply this theorem to our case, consider the one-sided sine-process on the space
of sequences $c_i\in\{0,1\}^{\mathbb{Z}^+} = X$. Set $c_i'=c_i-c_i c_{i+1}$ and
let $P$ be the corresponding probability distribution on the space of sequences
$c_i'$. Now let $T$ be the shift to the left, $(Tx)_i=x_{i+1}$, $f = c_0' -
\mathbb{E}(c_0')$ and let $\mathcal{F}_0$ be the standard $\sigma$-algebra
generated by cylinder sets.
$\mathbb{E}(\hat{T}\hat{T}^*\phi|\mathcal{F}_1)=\mathbb{E}(\phi|\mathcal{F}_1)$
because, if a function $g$ does not depend on the value $c_0'$, $\hat{T}^*g$ is just
the right shift of $g$.  Note that condition 1. is satisfied because
$$
\left|Cov(c_i',c_j')\right|\leq \frac{C}{(i-j)^2}
$$
for some $C>0$. To check the second condition note that for every cylinder set
the sum of considered series has a finite expectation because
$$
\left|Cov(c_{i+\vec{x}}, c_{j+\vec{y}})\right| \leq
\frac{C(\vec{x},\vec{y})}{(i-j)^2}
$$
where $\vec{x},\vec{y}$ are arbitrary finite subsets of $\mathbb{Z}^{+}$.

This immediately yields the central limit theorem for the Plancherel measure in the
regime of Proposition \ref{lv}:
\begin{proposition}
\label{clt}
 Assume that $\lim\limits_{n\to\infty}x_n = 1$ and $\forall\varepsilon > 0
\,|y_n-x_n|=o(n^{\varepsilon})$. Then the following holds:
$$
  \frac{\sum_{i \in [x_n,y_n]} c_i'(\lambda) - \mathbb{E}_{\Pl^{(n)}} (\sum_{i
\in [x_n,y_n]} c_i'(\lambda))}{\sqrt{Var_{\Pl^{(n)}}(\sum_{i \in [x_n,y_n]}
c_i'(\lambda))}} \xrightarrow[d]{} \mathcal{N}(0,1).
$$
\end{proposition}
Indeed, moments of the random variable $\frac{1}{\sqrt{n}}\sum_{i=0}^{n-1}\hat{T}^i
f$, in the notation introduced above, converges to the moments of
$\mathcal{N}(0,1)$.  Poissonized moments of $\sum{c_i'}$ converge to the moments
of the normal distribution from the Debye asymptotics (\ref{Debye}). Depoissonizing
inequalities in the Debye asymptotics one comes back to the Plancherel measure,
similarly to the depoissonization of \ref{moments}.
More precisely,  from the Debye asymptotics of the Bessel function one gets (see the proof
of Lemma 3.5 in \cite{key-2}):
$$
|J(x,y,\theta^2)-K_{sin}(x-y,\phi_x)|=O\left(\frac{|x-y|}{\sqrt{n}}
\right)\exp(\poissdev)
$$
where $\phi_x=\arccos(x/2\sqrt{n})$ is bounded away from the ends of the
interval $(0,\pi)$. Denoting
$$
X_n (\lambda) = \frac{\sum_{i \in [x_n,y_n]} c_i'(\lambda)}{\sqrt{Var_P(\sum_{i \in
[x_n,y_n]} c_i'(\lambda))}}
$$ 
we get
$$
|\mathbb{E}_{\theta^2}(X_n-\mathbb{E}_{\theta^2}(X_n))^k -
\mathbb{E}_P(X_n-\mathbb{E}_P(X_n))^k| = o(1)\exp(\poissdev) 
$$
Note that here we heavily depend on the fact that there are $O(n^k)$ summands of
the form $\prod_{i}K(i,\sigma(i))$, where $\sigma$ is some permutation, in the
expression for the $k^{th}$ moment of the
determinantal process with the kernel $K$.
Using the Debye asymptotics again, similarly to the depoissonization of the main
theorem and (\ref{moments}), we get
$$
|\mathbb{E}_{\theta^2}(X_n-\mathbb{E}_{\Pl^{(n)}}(X_n))^k - \mu_k| =
o(1)\exp(\poissdev),
$$ 
where $\mu_k$ is the $k^{th}$ moment of $\mathcal{N}(0,1)$. 
Depoissonizing and applying the Debye asymptotics again we get
$$
\lim\limits_{n\to\infty}\mathbb{E}_{\Pl^{(n)}}(X_n-\mathbb{E}_{\Pl^{(n)}}
(X_n))^k = \mu_k.
$$
This completes the proof of Proposition \ref{clt}.


\begin{thebibliography}{3}

\bibitem{key-5}Leonid Bogachev, Honggen Su, \emph{Central limit theorem
for random partitions under the Plancherel measure}, Doklady Mathematics, 2007,
v. 75, no. 3, pp. 381-384.

\bibitem{key-2}Alexei Borodin, Andrei Okounkov, Grigori Olshanski,
\emph{Asymptotics of Plancherel measures for symmetric groups, } J. Amer. Math.
Soc. 13 (2000), 481-515.

\bibitem{key-3}Patrik L. Ferrari, Alexei Borodin, \emph{Anisotropic
KPZ growth in 2+1 dimensions: fluctuations and covariance structure,
} J. Stat. Mech. (2009).

\bibitem{key-1}Alexander I. Bufetov, \emph{On the Vershik-Kerov
Conjecture Concerning the Shannon-Macmillan-Breiman Theorem for the
Plancherel Family of Measures on the Space of Young Diagrams}, arXiv:1001.4275

\bibitem{key-11}Costin, O., Lebowitz, J. L. \emph{Gaussian fluctuation in random
matrices.} Phys. Rev. Lett. 75,
69–72. (doi:10.1103/PhysRevLett.75.69).

\bibitem{key-9}Vladimir Ivanov, Grigori Olshanski, \emph{Kerov's central limit
theorem for the Plancherel measure on Young diagrams}, 	Symmetric Functions
2001: Surveys of Developments and Perspectives (NATO Science Series II.
Mathematics, Physics and Chemistry. Vol.74), Kluwer, 2002, pp. 93-151. 

\bibitem{key-8}Kurt Johansson, \emph{The longest increasing subsequence
in a random permutation and a unitary random matrix model}, Math. Res.
Letters, 5, 1998, 63\textendash{}82.

\bibitem{key-10}S. Kerov, \emph{Gaussian limit for the Plancherel measure of the
symmetric group}, Comptes
Rendus Acad. Sci. Paris, S\'{e}rie I 316 (1993), 303–308. 

\bibitem{key-7}Vershik, A. M.; Kerov, S. V. \emph{Asymptotic behaviour
of the Plancherel measure of the symmetric group and the limit form
of Young tableaux.} Dokl. Akad. Nauk SSSR 233 (1977), no. 6,
1024\textendash{}1027. 

\bibitem{key-14} Carlangelo Liverani, \emph{Central limit theorem for
deterministic systems}, International conference on dynamical systems
(Montevideo, 1995), 56-75.

\bibitem{key-6}Logan, B. F.; Shepp, L. A. \emph{A variational problem
for random Young tableaux.} Advances in Math. 26 (1977), no. 2,
206\textendash{}222.

\bibitem{key-4}Andrei Okounkov, \emph{Symmetric functions and
random partitions}, Symmetric functions 2001: surveys of developments
and perspectives, 223\textendash{}252, NATO Sci. Ser. II Math. Phys.
Chem., 74, Kluwer Acad. Publ., Dordrecht, 2002.

\bibitem{key-12}Alexander B. Soshnikov, \emph{Gaussian fluctuation for the
number of particles in Airy, Bessel, sine and other determinantal random point
fields}, J. Statist. Phys. 100 491-522. 

\bibitem{key-13}Alexander B. Soshnikov, \emph{Determinantal random point fields}
(Russian),
 Math. Surveys, 55, No. 5, 923-975, (2000).

\end{thebibliography}
\end{document}